\documentclass[12pt]{article}
\usepackage{epic,eepic,epsf,epsfig}
\usepackage{amsmath}
\usepackage{amssymb}
\usepackage{amsbsy}

\usepackage{amsfonts}
\newtheorem{theorem}{Theorem}[section]%
\newtheorem{lemma}[theorem]{Lemma}%
\newtheorem{cor}[theorem]{Corollary}%

\def\f{\noindent}

\def\mod{\hbox{\rm mod}\;}

\newcommand{\qed}{\mbox{\raisebox{0.7ex}{\fbox{}}} \vspace{4truemm}}
\def\demo{\f {\bf Proof.}\hskip10pt}

\begin{document}

\baselineskip 16pt

\title{ \vspace{-1.2cm}
On finite groups in which some maximal invariant subgroups have indices a prime or the square of a prime
\thanks{\scriptsize This research was supported in part by Shandong Provincial Natural Science Foundation, China (ZR2017MA022 and ZR2020MA044)
and NSFC (11761079).
\newline
 \hspace*{0.5cm} \scriptsize $^{\ast\ast}$Corresponding
  author.\newline
       \hspace*{0.5cm} \scriptsize{E-mail addresses:}
       shijt2005@163.com\,(J. Shi),\,m18863093906@163.com\,(Y. Tian).}}

\author{Jiangtao Shi\,$^{\ast\ast}$,\,Yunfeng Tian\\
\\
{\small School of Mathematics and Information Sciences, Yantai University, Yantai 264005, China}}

\date{ }

\maketitle \vspace{-.8cm}

\begin{abstract}
Let $A$ and $G$ be finite groups such that $A$ acts coprimely on $G$ by automorphisms, we first prove some results on the solvability of finite groups in which some maximal $A$-invariant subgroups have indices a prime or the square of a prime. Our results generalize Hall's theorem and some other known results. Moreover, we obtain a complete characterization of finite groups in which every non-nilpotent maximal $A$-invariant subgroup that contains the normalizer of some $A$-invariant Sylow subgroup has index a prime.

\medskip
\f {\bf Keywords:} maximal $A$-invariant subgroup; normalizer; non-nilpotent maximal $A$-invariant subgroup; solvable; normal $p$-complement\\
{\bf MSC(2020):} 20D10; 20D20
\end{abstract}

\section{Introduction}

In this paper all groups are assumed to be finite. It is known that Hall's theorem [4, Theorem 9.4.1] indicates that if every maximal subgroup of a group $G$ has index a prime or the square of a prime, then
$G$ is solvable. As a generalization of [4, Theorem 9.4.1], Shao and Beltr$\rm\acute{a}$n [5, Theorem A] gave the following result.

\begin{theorem} {\rm[5, Theorem A]}\ \ Let $G$ and $A$ be groups of coprime orders and assume that $A$ acts on $G$ by automorphisms. If the index of every maximal $A$-invariant
subgroup of $G$ is prime or the square of a prime, then $G$ is solvable.
\end{theorem}

In this paper, as a further generalization of [4, Theorem 9.4.1] and [5, Theorem A], we first obtain the following result whose proof is given in Section~3.

\begin{theorem}\ \ \label{th3} Let $A$ and $G$ be groups such that $A$ acts coprimely on $G$ by automorphisms. If
every non-nilpotent maximal $A$-invariant subgroup of $G$ that contains the normalizer of some $A$-invariant Sylow subgroup has index a prime or the square of a prime, then $G$ is solvable.
\end{theorem}

The following two corollaries are direct consequences of Theorem~1.2.

\begin{cor}\ \ \label{th2} Let $A$ and $G$ be groups such that $A$ acts coprimely on $G$ by automorphisms. If
every non-nilpotent maximal $A$-invariant subgroup of $G$ has index a prime or the square of a prime, then $G$ is solvable.
\end{cor}

\begin{cor}\ \ \label{th1} Let $A$ and $G$ be groups such that $A$ acts coprimely on $G$ by automorphisms. If
every maximal $A$-invariant subgroup of $G$ that contains the normalizer of some $A$-invariant Sylow subgroup has index a prime or the square of a prime, then $G$ is solvable.
\end{cor}

It is known that Huppert's theorem states that if every maximal subgroup of a group $G$ has
index a prime then $G$ is supersolvable, see [2, VI, Theorem 9.2]. As a generalization of Huppert's theorem, considering the indices of non-nilpotent maximal subgroups of a group,
Shi [6, Theorem 3] showed that a group in which every non-nilpotent maximal subgroup has prime index has a
Sylow tower. Beltr$\rm\acute{a}$n and Shao [1, Theorem C] and Li et al [3, Theorem 1.4] generalized [6, Theorem 3] as follows respectively.

\begin{theorem} {\rm[1, Theorem C]} \ \ Let $G$ and $A$ be groups of coprime orders and assume that $A$ acts on $G$ by automorphisms. If every non-nilpotent maximal $A$-invariant
subgroup of $G$ has prime index, then $G$ is solvable.
\end{theorem}

\begin{theorem} {\rm[3, Theorem 1.4]}\ \ Suppose that $A$ acts on $G$ via automorphisms and that $(|A|,|G|)=1$. If every non-nilpotent maximal $A$-invariant subgroup of $G$ has prime index, then
$G$ is a Sylow tower group.
\end{theorem}

As a further generalization of [6, Theorem 3] and [3, Theorem 1.4], we obtain the following result whose proof is given in Section~4.

\begin{theorem}\ \ \label{th4} Let $A$ and $G$ be groups such that $A$ acts coprimely on $G$ by automorphisms. If every non-nilpotent maximal $A$-invariant subgroup that contains the normalizer of some $A$-invariant Sylow subgroup has index a prime, then

$(1)$ $G$ has a Sylow tower;

$(2)$ Suppose that $p$ is the largest prime divisor of $|G|$ and $P$ is a $A$-invariant Sylow $p$-subgroup of $G$, then either $P$ is normal in $G$ or $G$
has a normal $p$-complement.
\end{theorem}

\section{Two Necessary Lemmas}

\begin{lemma} \ \ \label{l01} Let $p$ be the largest prime divisor of the order of a group $G$ and
$P\in{\rm Syl}_p(G)$, then either $P\trianglelefteq G$ or the maximal subgroup of $G$ that contains $N_G(P)$ has composite index.
\end{lemma}

\demo Suppose that $P$ is not normal in $G$. Let $M$ be maximal in $G$ satisfying $N_G(P)\leq M$ and $|G:M|=m$. In the following we will show that $m$ is a composite number.

Otherwise, if $m$ is a prime, then $m<p$ by the hypothesis. Since $|G:M|=m$, one has that $G/M_G$ is isomorphic to a subgroup of the symmetric group $S_m$,
where $M_G=\bigcap_{g\in G}M^g$ is the largest normal subgroup of $G$ that is contained in $M$. It is clear that $p\nmid|S_m|$, which implies that $p\nmid|G/M_G|$. It follows that $P\in{\rm Syl}_p(M_G)$.
By Frattini argument, one has $G=M_GN_G(P)\leq M_GM=M$, a contradiction. Therefore, $m$ is a composite number. \hfill\qed

\begin{lemma} {\rm[7, Lemma 9]}\ \ \label{l02} Let $H$ be a nilpotent Hall-subgroup of a group $G$ which is not a Sylow subgroup
of $G$. If for each prime divisor $p$ of $|H|$, assume $P\in{\rm Syl}_p(H)$ we always have $N_G(P)=H$, then there exists a normal
subgroup $K$ of $G$ such that $G=KH$ and $K\cap H=1$.
\end{lemma}

\section{Proof of Theorem~1.2}\label{s3}

\demo Let $G$ be a counterexample of minimal order. For any $A$-invariant normal subgroup $K$ of $G$, it is easy to see that the hypothesis of theorem holds for the quotient group $G/K$.
By the minimality of $G$, one has that $G/K$ is solvable. It follows that $K$ is non-solvable since $G$ is non-solvable.

Assume that $N$ is a minimal $A$-invariant normal subgroup of $G$. Let $p$ be the largest prime divisor of $|N|$, then $p$ is odd. Since $A$ acts coprimely on $G$ by automorphisms, take $P$ as a $A$-invariant Sylow $p$-subgroup of $G$. One has that $P_0=P\cap N$ is a $A$-invariant
Sylow $p$-subgroup of $N$. Since $G$ has no solvable $A$-invariant normal subgroup, one has $P_0\ntrianglelefteq G$. Then $N_G(P)\leq N_G(P_0)<G$. Note that $N_G(P_0)$ is $A$-invariant, assume that $M$
is a maximal $A$-invariant subgroup of $G$ such that $N_G(P)\leq N_G(P_0)\leq M$. By Frattini argument, one has $G=N_G(P_0)N=MN$.

For the maximal $A$-invariant subgroup $M$, we divide our analyses into two cases.

Case (1): Assume that $M$ is non-nilpotent, then $|G:M|=q$ or $r^2$ for some primes $q$ and $r$ by the hypothesis. Since $G=N_G(P_0)N=MN$, one has $|G:M|=|MN:M|=|N:N\cap M|=q$ or $r^2$. Note that $N\cap M\geq N\cap N_G(P_0)=N_N(P_0)$, so $q<p$ and $r<p$.

For the case when $|N:N\cap M|=q$, then $N\cap M$ is a maximal subgroup of $N$ that contains $N_N(P_0)$. It follows that $P_0\trianglelefteq N$ by Lemma~2.1, which implies that $P_0\trianglelefteq G$, a contradiction.

For another case when $|N:N\cap M|=r^2$. Note that $|N:N\cap M|=\frac{|N:N_N(P_0)|}{|N\cap M:N_N(P_0)|}$, $|N:N_N(P_0)|\equiv 1(\,\mod p)$ and $|N\cap M:N_N(P_0)|=|N\cap M:N_{(N\cap M)}(P_0)|\equiv1 (\,\mod p)$ by Sylow's theorem. Then we can easily get $p\mid r^2-1$. Since $r<p$, one has $p=r+1$, which implies that $r=2$ and $p=3$. Then $|N:N\cap M|=4$. It follows that $N/(N\cap M)_N$ is isomorphic to a subgroup of the symmetric subgroup $S_4$, where $(N\cap M)_N=\bigcap_{x\in N}(N\cap M)^x$ is the largest normal subgroup of $N$ that is contained in $N\cap M$. Note that $S_4$ is solvable, then $N'<N$. Since $N'\,{\rm char}\,N$ and $N$ is a $A$-invariant normal subgroup of $G$, one has that $N'$ is also a $A$-invariant normal subgroup of $G$. By the minimality of $N$, one has $N'=1$. It follows that $N$ is abelian, also a contradiction.

Case (2): Assume that $M$ is nilpotent. Since $G$ has no solvable $A$-invariant normal subgroup, $M$ must be a Hall-subgroup of $G$. Note that $p\mid|M|$ and $p$ is an odd prime, $M$ cannot be a Sylow subgroup of $G$ by
{\rm[1, Theorem B]}. Moreover, let $t$ be any prime divisor of $|M|$ and $T$ be a $A$-invariant Sylow $t$-subgroup of $M$, it is clear that $N_G(T)$ is $A$-invariant and $M\leq N_G(T)<G$.
By the maximality of $M$, one has $N_G(T)=M$. Therefore, there exists a normal Hall-subgroup $H$ of $G$ such that $G=M\ltimes H$ by Lemma~2.2. Since $M$ is a maximal $A$-invariant subgroup of $G$,
one has that $H$ is a minimal $A$-invariant normal subgroup of $G$. Note that $p\mid|N|$ but $p\nmid|H|$, then $N\neq H$. One has $H\cap N=1$. By the minimality of $G$, one has that both $G/H$ and $G/N$
are solvable. It follows that $G\cong G/(H\cap N)\lesssim (G/H)\times (G/N)$ and then $G$ is solvable, also a contradiction.

So the counterexample of minimal order does not exist and then $G$ is solvable.\hfill\qed

\section{Proof of Theorem~1.7}\label{s4}

\demo In the following we divide our discussions into two parts.

{\bf Part I.} We will prove that $G$ has a normal Sylow subgroup.

Let $G$ be a counterexample of minimal order.

Suppose that $Q_0$ is a minimal $A$-invariant normal subgroup of $G$. Since $G$ is solvable by Theorem~1.2, one has that $Q_0$ is a solvable characteristic simple group by the minimality of $Q_0$. Then $Q_0$ is an elementary group of prime-power order. Note that the hypothesis of the theorem also holds for $G/Q_0$, one has that $G/Q_0$
has a normal Sylow subgroup $SQ_0/Q_0$, where $S$ is a $A$-invariant Sylow subgroup of $G$. It follows that $SQ_0\trianglelefteq G$. By Frattini argument, one has $G=SQ_0N_G(S)=Q_0N_G(S)$.

If all maximal $A$-invariant subgroups of $G$ contain $Q_0$. Since $N_G(S)<G$ by the hypothesis, assume that $D$ is a maximal $A$-invariant subgroup of $G$ such that $N_G(S)\leq D$,
one has $G=Q_0N_G(S)=D$, a contradiction. Thus there exists a maximal $A$-invariant subgroup $M$ of $G$ such that $Q_0\nleq M$. It follows
that $G=Q_0M$.

Since $Q_0$ is an elementary abelian group, one has $Q_0\cap M\trianglelefteq Q_0$ and $Q_0\cap M\trianglelefteq M$, which implies that $Q_0\cap M\trianglelefteq Q_0M=G$.
Note that $Q_0\cap M<Q_0$, then $Q_0\cap M=1$ by the minimality of $Q_0$. So $G=Q_0\rtimes M$. Suppose $|Q_0|=q^\alpha$, where $\alpha\geq 1$.

For $M$ we divide our analyses into two cases.

For the case when $M$ is nilpotent, let $Q_1$ be a $A$-invariant Sylow $q$-subgroup of $M$. Then $Q_0\times Q_1\in{\rm Syl}_q(G)$ and $Q_0\times Q_1\trianglelefteq G$, a contradiction.

For the case when $M$ is non-nilpotent. Considering the quotient group $G/Q_0$. By the minimality of $G$, one has that $G/Q_0$ has a normal Sylow $r$-subgroup $RQ_0/Q_0$,
where $R$ is a $A$-invariant Sylow $r$-subgroup of $G$. Note that $G=Q_0\rtimes M$, without loss of generality, let $R$ be a $A$-invariant Sylow $r$-subgroup of $M$. Then $R=R(Q_0\cap M)=RQ_0\cap M\trianglelefteq M$. It follows
that $M\leq N_G(R)$. As $R\ntrianglelefteq G$, one has $M=N_G(R)$. Since $M$ is non-nilpotent, $|G:M|$ is a prime by the hypothesis. It follows that $|Q_0|=q$. Suppose $Q_0=Z_q$.

Let $p$ be the largest prime divisor of $|G|$ and let $P$ be a $A$-invariant Sylow $p$-subgroup of $G$. By our assumption, $P\ntrianglelefteq G$. Then each maximal $A$-invariant subgroup $K$ of $G$ that contains
$N_G(P)$ has composite index by Lemma~2.1. One has that $K$ is  nilpotent by the hypothesis.

If $Z_q\nleq K$, then $G=Z_q\rtimes K$. Arguing as above when $M$ is nilpotent, we can get a contradiction.

If $Z_q\leq K$, note that $K$ is nilpotent and $Z_q\trianglelefteq K$, one has $K\leq C_G(Z_q)$.

$(a)$ Suppose $K=C_G(Z_q)$. Since $N_G(Z_q)=G$, one has $K\trianglelefteq G$.
Take $T$ being a $A$-invariant Sylow $t$-subgroup of $K$, where $t\neq q$, then $T$ is a normal Sylow subgroup of $G$, a contradiction.

$(b)$ Suppose $K<C_G(Z_q)$. Since $K$ is a maximal $A$-invariant subgroup of $G$, one has
$C_G(Z_q)=N_G(Z_q)=G$, which implies that $Z_q\leq Z(G)$. Therefore, $G=Z_q\times M$. Note that $R\trianglelefteq M$, then $R$ is a normal Sylow subgroup of $G$,
also a contradiction.

All above analyses imply that the counterexample of minimal order does not exist and so $G$ has a normal Sylow subgroup.

Moreover, since the hypothesis of the theorem also holds for the quotient group $G/N$ for any $A$-invariant normal subgroup $N$ of $G$, one can easily get that $G$ has a Sylow tower.

{\bf Part II.} Let $p$ be the largest prime divisor of $|G|$, we will prove that either the $A$-invariant Sylow $p$-subgroup of $G$ is normal or $G$ has a normal $p$-complement.

Let $G$ be a counterexample of minimal order. Then the $A$-invariant Sylow $p$-subgroup of $G$ is not normal.

By Part I, let $Q$ be a normal Sylow $q$-subgroup of $G$, where $q\neq p$. Since $G$ is solvable by Theorem~1.2, there exists a $A$-invariant subgroup $L$ of $G$ such that $G=Q\rtimes L$ by Schur-Zassenhaus theorem.
Let $P$ be a $A$-invariant Sylow $p$-subgroup of $L$, then $P$ is also a $A$-invariant Sylow $p$-subgroup of $G$. For the quotient group $G/Q$, $p$ is also the largest prime divisor of $|G/Q|$ and $PQ/Q$ is a $A$-invariant Sylow $p$-subgroup of $G/Q$.
By the minimality of $G$, one has either $PQ/Q\trianglelefteq G/Q$ or $G/Q$ has a normal $p$-complement $U/Q$.

$(i)$ If $PQ/Q\trianglelefteq G/Q$, then $PQ\trianglelefteq G$. It follows that $PQ\cap L\trianglelefteq L$ and then $P=P(Q\cap L)=PQ\cap L\trianglelefteq L$.
It implies that $L\leq N_G(P)<G$. By the hypothesis, one has that each maximal $A$-invariant subgroup of $G$ that contains $N_G(P)$ is nilpotent. Then $L$ is also nilpotent.
Assume $L=P\times L_1$, then $P$ has a normal $p$-complement $Q\times L_1$, a contradiction.

$(ii)$ If $G/Q$ has a normal $p$-complement $U/Q$, then $G/Q=(PQ/Q)\ltimes (U/Q)$, where $p\nmid |U|$. One has $G=PQU=PU=P\ltimes U$, which implies that $U$ is a
normal $p$-complement of $G$, also contradiction.

Hence the counterexample of minimal order does not exist and so either the $A$-invariant Sylow $p$-subgroup of $G$ is normal or $G$ has a normal $p$-complement.\hfill\qed

\end{document}